\documentclass{article}
\usepackage[english]{babel}
\usepackage{amsmath,amsthm,amssymb,enumerate,bbm}
\usepackage{geometry}
\catcode`\<=\active \def<{
\fontencoding{T1}\selectfont\symbol{60}\fontencoding{\encodingdefault}}
\catcode`\>=\active \def>{
\fontencoding{T1}\selectfont\symbol{62}\fontencoding{\encodingdefault}}
\catcode`\|=\active \def|{
\fontencoding{T1}\selectfont\symbol{124}\fontencoding{\encodingdefault}}
\newcommand{\assign}{:=}
\newcommand{\dueto}[1]{\textup{\textbf{(#1) }}}
\newcommand{\longhookrightarrow}{{\lhook\joinrel\relbar\joinrel\rightarrow}}
\newcommand{\mathd}{\mathrm{d}}
\newcommand{\nocomma}{}

\newcommand{\tmdate}[1]{\today}
\newcommand{\tmmathbf}[1]{\ensuremath{\boldsymbol{#1}}}
\newcommand{\tmop}[1]{\ensuremath{\operatorname{#1}}}
\newcommand{\tmtextit}[1]{{\itshape{#1}}}
\newenvironment{enumeratealpha}{\begin{enumerate}[a{\textup{)}}] }{\end{enumerate}}
\newenvironment{enumerateroman}{\begin{enumerate}[i.] }{\end{enumerate}}
\newtheorem{corollary}{Corollary}
\newtheorem{lemma}{Lemma}
\newtheorem{remark}{Remark}
\newtheorem{theorem}{Theorem}

\begin{document}

\title{On Stationary Solutions of the 2D Doi-Onsager Model}

\author{
  Mohammad Ali Niksirat, Xinwei Yu
}

\maketitle

\begin{abstract}
  We study the 2D Doi--Onsager models with general potential kernel, with
  special emphasis on the classical Onsager kernel. Through application of
  topological methods from nonlinear functional analysis, in particular the
  Leray--Schauder degree theory, we obtain the uniqueness of the trivial
  solution for low temperatures as well as the local bifurcation structure of
  the solutions.
\end{abstract}

\section{Introduction}

In 1949, Lars Onsager proposed a mathematical model for the phase transition
of equilibria of dilute colloidal solutions of rod-like molecules between the
isotropic and nematic phases ({\cite{Onsager1949}}). As the fluid in both
phases is homogeneous, that is the locations of the molecules do not matter,
Onsager's theory focuses on a probability density function $f ( \tmmathbf{r}
)$ over the unit sphere which models distribution of the directions of the
rods. Although the original modeling is carried out in $\mathbbm{R}^{3}$, the
mathematical formulation can be generalized to $\mathbbm{R}^{d}$ for any
dimension $d \geqslant 2$ in a straightforward manner. In the following we
present this generalized version.

Denote by $S^{d-1}$ the unit sphere in $\mathbbm{R}^{d}$. Let $f (
\tmmathbf{r} ) :S^{d-1} \mapsto [ 0, \infty )$ be the probability density
characterizing the directions of the rods, that is
\begin{equation}
  P \left( \text{the rod is along } \tmmathbf{r} \in A \subseteq S^{d-1}
  \right) = \int_{A} f ( \tmmathbf{r} ) \mathd \sigma ( \tmmathbf{r} )
\end{equation}
where we denote by $\sigma ( \tmmathbf{r} )$ the volume element on $S^{d-1}$.
As we are modeling ``rod-like'' molecules with no distinction between the two
ends, we can further assume $f ( \tmmathbf{r} ) =f ( -\tmmathbf{r} )$.
Consequently the constraints on $f ( \tmmathbf{r} )$ are
\begin{equation}
  f ( \tmmathbf{r} ) \geqslant 0, \hspace{2em} f ( \tmmathbf{r} ) =f (
  -\tmmathbf{r} ) , \hspace{2em} \int_{S^{d-1}} f ( \tmmathbf{r} ) \mathd
  \sigma ( \tmmathbf{r} ) =1. \label{eq:201302151142}
\end{equation}

The equilibrium distributions correspond to the critical points of the
following functional:
\begin{equation}
  E ( f ) := \int_{S^{d-1}} ( \log  f ( \tmmathbf{r} ) ) f ( \tmmathbf{r} )
  \mathd \sigma ( \tmmathbf{r} ) + \frac{1}{2}  \int_{S^{d-1}} ( U ( f ) (
  \tmmathbf{r} ) ) f ( \tmmathbf{r} ) \mathd \sigma ( \tmmathbf{r} )
  \label{eq:201302151138}
\end{equation}
which is derived in {\cite{Onsager1949}} as the second order approximation of
the free energy -- neglecting interactions between three and more molecules.

The interaction potential $U ( f )$ in (\ref{eq:201302151138}) is given by
\begin{equation}
  U ( f ) ( \tmmathbf{r} ) \assign \lambda \int_{S^{d-1}} K (
  \tmmathbf{r},\tmmathbf{r}' ) f ( \tmmathbf{r}' ) \mathd \sigma (
  \tmmathbf{r}' ) \label{eq:201302151139}
\end{equation}
where the parameter $\lambda >0$ can be interpreted as either the
concentration of the particles in the carrier fluid or the inverse of the
absolute temperature. The interaction kernel $K ( \tmmathbf{r},\tmmathbf{r}'
)$ inherits the following symmetry properties.
\begin{equation}
  K ( -\tmmathbf{r},\tmmathbf{r}' ) =K ( \tmmathbf{r},\tmmathbf{r}' ) ;
  \hspace{1em} K ( \tmmathbf{r},\tmmathbf{r}' ) =K ( \tmmathbf{r}'
  ,\tmmathbf{r} ) ; \hspace{1em} K ( \tmmathbf{r},\tmmathbf{r}' ) =K (
  T\tmmathbf{r},T\tmmathbf{r}' ) \hspace{2em} \forall T \in O ( 3 ) .
  \label{eq:201302151140}
\end{equation}
Note that when $d=2$ we can use the natural parametrization of $S^{1}$ by the
angle $\theta \in [ 0,2 \pi )$ and rewrite any kernel satisfying
(\ref{eq:201302151140}) as a convolution kernel $K ( \theta - \theta' )$ for
some even function $K$ satisfying $K ( \theta + \pi ) =K ( \theta )$. This
reduces the right hand side of (\ref{eq:201302151139}) to a convolution
\begin{equation}
  U ( f ) ( \theta ) \assign \lambda \int_{0}^{2 \pi} K ( \theta - \theta' ) f
  ( \theta' ) \mathd \theta' . \label{eq:201304021253}
\end{equation}

The Euler--Lagrange equation for the system
(\ref{eq:201302151142})--(\ref{eq:201302151139}) can be easily written down as
\begin{equation}
  f ( \tmmathbf{r} ) = \frac{e^{-U ( f ) ( \tmmathbf{r} )}}{\int_{S^{d-1}}
  e^{-U ( f ) ( \tmmathbf{r} )} \mathd \sigma} , \hspace{2em} f ( \tmmathbf{r}
  ) =f ( -\tmmathbf{r} ) . \label{eq:201302151143}
\end{equation}
A moment's inspection reveals that $f ( \tmmathbf{r} ) \equiv \frac{1}{|
S^{n-1} |}$ is always a solution. This constant solution corresponds to the
uniform distribution of rod directions and therefore models the ``isotropic''
or ``un-ordered'' phase where all directions are equally likely to be taken by
the molecules. On the other hand, it has been observed since 1888
({\cite{Reinitzer1888}}) that as the temperature $\lambda^{-1}$ decreases, the
fluid may go through one or more phase transitions resulting in some order of
the directions taken by the molecules. Such phase transition to the so-called
nematic phases can be modeled by the bifurcation of the constant solution to
non-constant solutions of (\ref{eq:201302151142})--(\ref{eq:201302151139}), or
equivalently of (\ref{eq:201302151143}).

The original kernel proposed by Onsager is
\begin{equation}
  K ( \tmmathbf{r},\tmmathbf{r}' ) = | \sin   \theta | \hspace{1em} \left( = |
  \tmmathbf{r} \times \tmmathbf{r}' | \text{ when } d=3 \right)
  \label{eq:201302151141}
\end{equation}
where $\theta$ is the angle between the unit vectors $\tmmathbf{r}$ and
$\tmmathbf{r}'$. For (\ref{eq:201302151143}) with this kernel, Onsager showed
through asymptotic expansion in {\cite{Onsager1949}} that when $\lambda$ is
large enough, bifurcation to non-constant solutions occur.

More quantitative analysis of the system (\ref{eq:201302151143}) with Onsager
kernel turned out to be difficult. On the other hand there are kernels
capturing the qualitative behavior of the solution while at the same time are
more friendly to mathematical analysis. One such kernel, due to Maier and
Saupe {\cite{Maier1958}}, reads
\begin{equation}
  K ( \tmmathbf{r},\tmmathbf{r}' ) = \cos^{2} \theta - \frac{1}{3} = (
  \tmmathbf{r} \cdot \tmmathbf{r}' )^{2} - \frac{1}{3} .
  \label{eq:201302151150}
\end{equation}
The Maier--Saupe kernel is often simply written as $( \tmmathbf{r} \cdot
\tmmathbf{r}' )^{2}$ as (\ref{eq:201302151143}) remains the same if we discard
the constant $- \frac{1}{3}$.

The major difference between (\ref{eq:201302151143}) with Maier--Saupe
potential (\ref{eq:201302151150}) and that with the Onsager potential
(\ref{eq:201302151141}) is that for the former the potential $U ( f )$, given
by (\ref{eq:201302151139}), resides in a finite dimensional space, thus
reducing the infinite dimensional problem (\ref{eq:201302151143}) to a finite
dimensional nonlinear system of equations. This reduced system, still highly
nontrivial, is nevertheless more tractable than the original system. As a
consequence, (\ref{eq:201302151143}) with Maier--Saupe potential has been well
understood through brilliant work of many researchers (see
{\cite{Constantin2004}}, {\cite{Fatkullin2005a}}, {\cite{Liu2005}},
{\cite{Zhou2005}}, {\cite{Liu2007}}, {\cite{Zhou2007}} for the case $d=3$,
{\cite{Constantin2005}}, {\cite{Fatkullin2005}}, {\cite{Luo2005}} for the case
$d=2$, and {\cite{Wang2008a}} for the general $d$-dimensional case.) Inspired
by these works, (\ref{eq:201302151143}) with other kernels enjoying similar
``dimension-reduction'' property has also be analyzed, see e.g.
{\cite{Chen2010}}. \

With the Maier--Saupe model (\ref{eq:201302151143}), (\ref{eq:201302151150})
understood, interest in the original Onsager model
(\ref{eq:201302151143})--(\ref{eq:201302151141}) was resurrected. Much
progress has been made in the past few years in the case $d=2$. In
{\cite{Chen2010}}, the axisymmetry of all possible solutions is proved, that
is, for any solution $f ( \theta )$ to
(\ref{eq:201302151143})--(\ref{eq:201302151141}), there is $\theta_{0}$ such
that $f ( \theta_{0} - \theta ) =f ( \theta_{0} + \theta )$. It is also proved
in {\cite{Chen2010}} that for appropriate $\lambda$, there are solutions of
arbitrary periodicity. In {\cite{Wang2008}} the authors rewrite
(\ref{eq:201302151143}) into an infinite system of nonlinear equations for the
Fourier coefficients of $f ( \theta )$ and calculated numerically the first
few bifurcations. More recently, in {\cite{Lucia2010}} the authors study the
case $d=2$ through cutting-off
(\ref{eq:201302151143})--(\ref{eq:201302151141}) to a finite dimensional
system of nonlinear equations, and obtain local bifurcation structure for this
finite dimensional approximation.

In this article, we try to gain more understanding of the original infinite
dimensional problem (\ref{eq:201302151143})--(\ref{eq:201302151141}) in the
case $d=2$:
\begin{equation}
  f ( \theta ) = \frac{e^{-U ( f ) ( \theta )}}{\int_{0}^{2 \pi} e^{-U ( f ) (
  \theta )} \mathd \theta} , \hspace{1em} f ( \theta ) =f ( \pi + \theta ) ,
  \hspace{1em} U ( f ) ( \theta ) = \lambda \int_{0}^{2 \pi} K ( \theta -
  \theta' ) f ( \theta' ) \mathd \theta' . \label{eq:201303271306}
\end{equation}
We show that most of the results obtained in {\cite{Lucia2010}} for the finite
dimensional truncated system of (\ref{eq:201303271306}) can be generalized to
the original infinite dimensional system (\ref{eq:201303271306}) itself. More
specifically, we have the following results.

Let $k_{m}$, $m=1,2,3, \ldots$ be defined through the Fourier expansion
\begin{equation}
  K ( \theta ) = \sum_{m=0}^{\infty} k_{m} \cos ( 2m \theta ) .
\end{equation}
\begin{itemize}
  \item {\dueto{Theorem \ref{thm:201304021429}}} The problem has a unique
  solution, which must be the constant solution, when $0< \lambda <
  \lambda_{0} \assign \left( \sum_{m=1}^{\infty} | k_{m} | \right)^{-1}$. This
  generalizes Proposition 3.1 b) in {\cite{Lucia2010}}.
  
  \item {\dueto{Theorem \ref{thm:201502121119}}}Two solutions bifurcate from
  the trivial solution at every $\lambda_{m} \assign - \frac{2}{k_{m}}$. The
  bifurcation is supercritical if $\frac{2k_{2m}}{k_{m}} <1$ and subcritical
  if $\frac{2k_{2m}}{k_{m}} >1$. Furthermore, in the former case the first
  pair of bifurcated solutions are stable and the other bifurcated solutions
  are unstable, while in the latter case all bifurcated solutions are
  unstable. This generalizes Proposition 4.4 and Corollary 4.5 in
  {\cite{Lucia2010}}.
\end{itemize}
Application of these results to the equation with Onsager's kernel leads to
the following conclusions.
\begin{itemize}
  \item The problem has a unique (trivial) solution when $0< \lambda <
  \frac{\pi}{2}$.
  
  \item Two solutions bifurcate from the trivial solution at $\lambda_{m} =
  \frac{( 4m^{2} -1 ) \pi}{2}$, $m=1,2,3, \ldots$. All bifurcations are
  supercritical.
  
  \item The pair of solutions bifurcating from $\lambda_{1} = \frac{3
  \pi}{2}$ is stable. All other bifurcated solutions are unstable. 
\end{itemize}
\begin{remark}
  \label{rem:201304021310}Our method applies in principle to the general cases
  $d \geqslant 3$ as well. However some technical difficulties arise and many
  new measures need to be taken. We will report our effort in this direction
  in a forthcoming paper.
\end{remark}

The remaining of the paper is organized as follows. In Section
\ref{sec:201303261407} we rewrite the problem (\ref{eq:201303271306}) into a
new formulation better-suited for the application of topological methods, and
carry out the calculation of the Jacobian matrix of the linearized operator.
In Section \ref{sec:201303261452} we prove that for all $0< \lambda <
\lambda_{0}$ the problem has a unique solution, which is trivial. In Section
\ref{sec:201303261453} we study the local bifurcation structure of the
problem. To improve the readability of the paper, statements of classical
results as well as some detailed calculations are delegated to Appendix
\ref{app:201303261454}.

\section{Preparations}\label{sec:201303261407}

\subsection{Re-formulation of the Problem}

Recall that we need to solve
\begin{equation}
  f ( \tmmathbf{r} ) = \frac{e^{-U ( f ) ( \tmmathbf{r} )}}{\int_{S^{d-1}}
  e^{-U ( f ) ( \tmmathbf{r} )} \mathd \sigma ( \tmmathbf{r} )} , \hspace{2em}
  f ( \tmmathbf{r} ) =f ( -\tmmathbf{r} ) . \label{eq:201303261410}
\end{equation}
with
\begin{equation}
  U ( f ) ( \tmmathbf{r} ) = \lambda \int_{S^{d-1}} K (
  \tmmathbf{r},\tmmathbf{r}' ) f ( \tmmathbf{r}' ) \mathd \sigma (
  \tmmathbf{r}' ) \label{eq:201303261411}
\end{equation}
Multiplying both sides of (\ref{eq:201303261410}) by $\lambda K (
\tmmathbf{r},\tmmathbf{r}' )$ and integrating over $S^{d-1}$, we cancel $f$
and reach an equation for the potential $U ( \tmmathbf{r} )$.
\begin{equation}
  U ( \tmmathbf{r} ) = \frac{\int_{S^{d-1}} \lambda K (
  \tmmathbf{r},\tmmathbf{r}' ) e^{-U ( \tmmathbf{r}' )} \mathd \sigma (
  \tmmathbf{r}' )}{\int_{S^{d-1}} e^{-U ( \tmmathbf{r} )} \mathd \sigma (
  \tmmathbf{r} )} , \hspace{2em} U ( \tmmathbf{r} ) =U ( -\tmmathbf{r} ) .
  \label{eq:201303140907}
\end{equation}
Note that once (\ref{eq:201303140907}) is solved, $f ( \tmmathbf{r} )$ can be
recovered from
\begin{equation}
  f ( \tmmathbf{r} ) = \frac{e^{-U ( \tmmathbf{r} )}}{\int_{S^{d-1}} e^{-U (
  \tmmathbf{r} )} \mathd \sigma ( \tmmathbf{r} )} . \label{eq:201304021327}
\end{equation}
Thus (\ref{eq:201303140907}) is equivalent to the original problem
(\ref{eq:201303261410})--(\ref{eq:201303261411}).

From now on we restrict ourselves to the specific case $d=2$. In this case we
can apply the natural parametrization of $S^{1}$ and write $K (
\tmmathbf{\tmmathbf{r},\tmmathbf{r}'} )$ as a convolution kernel $K ( \theta -
\theta' )$. This reduces (\ref{eq:201303140907}) to
\begin{equation}
  U ( \theta ) = \frac{\int_{0}^{2 \pi} \lambda K ( \theta - \theta' ) e^{-U (
  \theta' )} \mathd \theta'}{\int_{0}^{2 \pi} e^{-U ( \theta )} \mathd \theta}
  , \hspace{2em} U ( \theta ) =U ( \theta + \pi ) . \label{eq:201502111419}
\end{equation}
Now we define $\overline{K} \assign \frac{1}{2 \pi}  \int_{0}^{2 \pi} K (
\theta ) \mathd \theta$ and denote
\begin{equation}
  \tilde{K} ( \theta ) \assign K ( \theta ) - \overline{K} , \hspace{2em} V (
  \theta ) \assign U ( \theta ) - \lambda \overline{K} .
\end{equation}
It is easy to see that (\ref{eq:201502111419}) is equivalent to the following.
\begin{equation}
  V ( \theta ) = \lambda \Gamma ( V ) ( \theta ) \assign \frac{\lambda
  \int_{0}^{2 \pi} \tilde{K} ( \theta - \theta' ) e^{-V ( \theta' )} \mathd
  \theta'}{\int_{0}^{2 \pi} e^{-V ( \theta )} \mathd \theta} , \hspace{1em}
  \int_{0}^{2 \pi} V ( \theta ) \mathd \theta =0, \hspace{1em} V ( \theta ) =V
  ( \theta + \pi ) . \label{eq:201303261447}
\end{equation}
As the kernel has rotational invariance and the solution is axisymmetric
({\cite{Chen2010}}), we can further require $V ( \theta ) =V ( 2 \pi - \theta
)$. We also assume that $K ( \theta ) \in W^{1, \infty} ( [ 0,2 \pi ] )$. Note
that this assumption is satisfied by all the kernels proposed in the
literature. The natural function space we will be working in is
\begin{equation}
  H \assign \left\{ V ( \theta ) \in H^{1} ( [ 0,2 \pi ] ) ;V ( \theta ) =V (\theta + \pi )  a.e.; \int_{0}^{2 \pi} V ( \theta ) \mathd \theta =0;V (\theta ) =V ( 2 \pi - \theta )  a.e. \right\} . \label{eq:201303140924}
\end{equation}

To summarize, we will study the fixed-point problem
\begin{equation}
  V ( \theta ) = \lambda \Gamma ( V ) ( \theta ) , \hspace{2em} V ( \theta )
  \in H,
\end{equation}
where
\begin{equation}
  \Gamma ( V ) ( \theta ) \assign \frac{\int_{0}^{2 \pi} \tilde{K} ( \theta -
  \theta' ) e^{-V ( \theta' )} \mathd \theta'}{\int_{0}^{2 \pi} e^{-V ( \theta
  )} \mathd \theta} . \label{eq:201502111432}
\end{equation}
\subsection{The Jacobian $D \Gamma$}\label{subsec:201303271012}

We calculate the Jacobian matrix $( a_{m n} ) =A \assign D  \Gamma$.

Denote $\phi_{n} := \frac{1}{\sqrt{( 4n^{2} +1 ) \pi}} \cos ( 2n \theta )$
which with $n=1,2,3, \ldots$ form an orthonormal basis for $H$. Then standard
calculation gives
\begin{eqnarray}
  D \Gamma ( V ) ( U ) ( \theta ) & = & \frac{\left( \int_{0}^{2 \pi}
  \tilde{K} ( \theta - \theta' ) e^{-V ( \theta' )} \mathd \theta' \right) 
  \int_{0}^{2 \pi} U ( \theta ) e^{-V ( \theta )} \mathd \theta}{\left(
  \int_{0}^{2 \pi} e^{-V ( \theta )} \mathd \theta \right)^{2}} \nonumber\\
  &  & - \frac{\int_{0}^{2 \pi} \tilde{K} ( \theta - \theta' ) U ( \theta' )
  e^{-V ( \theta' )} \mathd \theta'}{\int_{0}^{2 \pi} e^{-V ( \theta )} \mathd
  \theta} .  \label{eq:201303271015}
\end{eqnarray}
If we define the probability measure
\begin{equation}
  \mathd \mu_{V} \assign \left( \int_{0}^{2 \pi} e^{-V ( \theta )} \mathd
  \theta \right)^{-1} e^{-V ( \theta )} \mathd \theta \label{eq:201303271124}
\end{equation}
then we can simplify (\ref{eq:201303271015}) to
\begin{equation}
  D \Gamma ( V ) ( U ) = \left[ \int_{0}^{2 \pi} \tilde{K} ( \theta - \theta'
  ) \mathd \mu_{V} ( \theta' ) \cdot \int_{0}^{2 \pi} U ( \theta' ) \mathd
  \mu_{V} ( \theta' ) - \int_{0}^{2 \pi} \tilde{K} ( \theta - \theta' ) U (
  \theta' ) \mathd \mu_{V} ( \theta' ) \right] \label{eq:201304231123}
\end{equation}
which leads to
\begin{eqnarray}
  A \assign ( a_{m n} ) & \assign & ( D \Gamma ( V ) ( \phi_{n} ) , \phi_{m}
  )_{H} \nonumber\\
  & = & \int_{0}^{2 \pi} ( D \Gamma ( V ) ( \phi_{n} ) ) \phi_{m} \mathd
  \theta + \int_{0}^{2 \pi} ( D \Gamma ( V ) ( \phi_{n} ) )' \phi_{m}' \mathd
  \theta \nonumber\\
  & = & \frac{k_{m} A_{m n}  ( 1+4mn )}{\sqrt{4m^{2} +1}  \sqrt{4n^{2} +1}} .
\end{eqnarray}
where $k_{m }$ are the Fourier coefficients of $\tilde{K} ( \theta )$
\begin{equation}
  \tilde{K} ( \theta ) \assign \sum_{m=1}^{\infty} k_{m} \cos ( 2m \theta )
\end{equation}
and
\begin{equation}
  A_{m n} \assign \int_{0}^{2 \pi} \cos ( 2m \theta ) \mathd \mu_{V} ( \theta
  ) \cdot \int_{0}^{2 \pi} \cos ( 2n \theta ) \mathd \mu_{V} ( \theta ) -
  \int_{0}^{2 \pi} \cos ( 2m \theta ) \cos ( 2n \theta ) \mathd \mu_{V} (
  \theta ) , \label{eq:201502121137}
\end{equation}

An important property of the matrix $A$ is $| a_{m n} | \leqslant | k_{m} |$.
To see this, we apply Lemma \ref{lem:201304181013} (see Appendix
\ref{subsec:201502101458}) to (\ref{eq:201502121137}) to conclude $| A_{m n} |
\leqslant 1$, from which the conclusion immediately follows.

\section{Uniqueness of the Trivial Solution}\label{sec:201303261452}

In this section we prove the following theorem.

\begin{theorem}
  \label{thm:201304021429}Assume $K \in W^{1, \infty} ( [ 0,2 \pi ] )$. Let
  $\lambda \sum_{m=1}^{\infty} | k_{m} | <1$. Then $V=0$, that is the only
  solution for (\ref{eq:201302151143})--(\ref{eq:201302151141}) is the
  constant solution. Here $k_{m}$ is the $m$-th coefficient of the Fourier
  expansion of $K ( \theta )$.
  \begin{equation}
    K ( \theta ) = \sum_{m=0}^{\infty} k_{m} \cos ( 2m \theta ) .
  \end{equation}
\end{theorem}

\begin{remark}
  This is a direct generalization of Proposition 3.1 b) of {\cite{Lucia2010}}
  to the infinite dimensional case. 
\end{remark}

The proof applies the classical Leray--Schauder theory to the fixed point
problem
\begin{equation}
  ( I- \lambda \Gamma ) ( V ) =0, \hspace{2em} V \in H \label{eq:201303140929}
\end{equation}
where $\Gamma$ is defined in (\ref{eq:201303261447}) and the space $H$ is
defined in (\ref{eq:201303140924}). To do this we need $H$ to be Hilbert and
$\Gamma$ to be compact, which are established by the following lemmas whose
proofs are delegated to Appendix \ref{app:201303270947}.

\begin{lemma}
  \label{lem:201304161059}$H$ is a Hilbert space. Furthermore $\Gamma :H
  \mapsto H$ if $K \in W^{1, \infty} ( [ 0,2 \pi ] )$.
\end{lemma}

\begin{lemma}
  \label{lem:201303270911} Assume $K ( \theta ) \in C ( [ 0,2 \pi ] )$. Then
  $\Gamma :H \mapsto H$ is compact.
\end{lemma}

\begin{remark}
  We emphasize that since $\Gamma$ is nonlinear, compactness here means (see
  e.g. {\cite{Nirenberg2001}})
  \begin{enumerateroman}
    \item $\Gamma$ is continuous;
    
    \item For every bounded closed $\Omega \subset H$, $\overline{\Gamma (
    \Omega )}$ is compact. 
  \end{enumerateroman}
\end{remark}

\noindent\textbf{Proof of Theorem \ref{thm:201304021429}.}
  As $W^{1, \infty} ( [ 0,2 \pi ] ) \hookrightarrow C ( [ 0,2 \pi ] )$, we can
  apply Lemma \ref{lem:201303270911} to conclude that $\Gamma$ is compact.
  
  The proof will now be carried out as follows. First we show the existence
  of a bounded open set $\Omega \subset H$ such that there is no solution
  outside $\Omega$. Next we show that the degree $\deg ( I- \lambda \Gamma ,
  \Omega ,0 ) =1$. Finally we prove that any possible solution to $( I-
  \lambda \Gamma ) ( V ) =0$ is isolated with index $1$. As in this case the
  degree is the sum of indices, we know that $0$ is the only solution.
  \begin{itemize}
    \item The existence of a bounded open set $\Omega \subset H$ such that $(
    I- \lambda \Gamma ) ( V ) =0$ has no solution outside $\Omega$. Let $R
    \assign \| K \|_{W^{1, \infty}} / \sum_{m=1}^{\infty} | k_{m} |$. Then it
    is easy to see that $\| \lambda \Gamma ( V ) \|_{H} \leqslant CR$ for all
    $\lambda$ satisfying the assumption of the theorem. Thus we can take
    $\Omega \assign B_{CR}$, the ball centered at the origin with radius $CR$.
    
    \item $\deg ( I- \lambda \Gamma , \Omega ,0 ) =1$.
    
    Introduce the homotopy $H ( t ) \assign I-t \lambda \Gamma$ with $t \in [
    0,1 ]$. We easily verify that $H ( t ) ( V ) =0$ has no solution on
    $\partial \Omega$ for all $t \in [ 0,1 ]$. Consequently
    \begin{equation}
      \deg ( I- \lambda \Gamma , \Omega ,0 ) = \deg ( H ( 1 ) , \Omega ,0 ) =
      \deg ( H ( 0 ) , \Omega ,0 ) = \deg ( I, \Omega ,0 ) =1.
    \end{equation}
    \item The solutions are isolated.
    
    The Fr{\'e}chet differentiability of $\Gamma$ can be verified through
    straightforward calculation, taking advantage of the embedding $H
    \longhookrightarrow L^{\infty} ( [ 0,2 \pi ] )$. The solutions are
    isolated if we can show that $I- \lambda D \Gamma$ is a homeomorphism. As
    $\lambda \Gamma$ is compact, so is the derivative $\lambda D \Gamma$.
    Applying standard Fredholm alternative (see e.g. {\cite{Ambrosetti1993}})
    we see that all we need to show is that $\ker ( I- \lambda D \Gamma ) = \{
    0 \}$.
    
    Take any $U \in \ker ( I- \lambda D \Gamma ) ( V )$. We have, following
    (\ref{eq:201304231123}) in Section \ref{subsec:201303271012},
    \begin{equation}
      U ( \theta ) = \lambda \left[ \int_{0}^{2 \pi} \tilde{K} ( \theta -
      \theta' ) \mathd \mu_{V} ( \theta' ) \cdot \int_{0}^{2 \pi} U ( \theta'
      ) \mathd \mu_{V} ( \theta' ) - \int_{0}^{2 \pi} \tilde{K} ( \theta -
      \theta' ) U ( \theta' ) \mathd \mu_{V} ( \theta' ) \right]
    \end{equation}
    where $\mathd \mu_{V}$ is as defined in (\ref{eq:201303271124}).
    
    Application of Lemma \ref{lem:201304181013} gives
    \begin{equation}
      | U ( \theta ) | \leqslant \lambda \| \tilde{K} \|_{L^{\infty}}  \| U
      \|_{L^{\infty}} . \hspace{2em} \forall \theta \in [ 0,2 \pi ]
    \end{equation}
    By assumption $\lambda \sum_{m=1}^{\infty} | k_{m} | <1$ which leads to
    $\lambda \| \tilde{K} \|_{L^{\infty}} <1$, consequently $U=0$.
    
    \item The index of any solution is $1$.
    
    Following the calculation in Section \ref{subsec:201303271012} we have $|
    a_{m n} | \leqslant | k_{m} |$ where $( a_{m n} )$ is the infinite
    dimensional matrix representation of $D  \Gamma$ with respect to the
    orthonormal basis $\left\{ \frac{1}{\sqrt{( 4n^{2} +1 ) \pi}} \cos ( 2n
    \theta ) \right\}_{n=1}^{\infty}$ of $H$. By assumption
    $\sum_{m=1}^{\infty} \lambda | k_{m} | <1$, therefore the eigenvalues of
    $I- \lambda D \Gamma$ are all bounded below by a positive constant.
    Consequently the index of the map $I- \lambda \Gamma$ is $1$ everywhere.
  \end{itemize}
  Thus we see that the desired conclusion holds when $\lambda \| \tilde{K}
  \|_{L^{\infty}} <1$ and $\sum_{m=1}^{\infty} \lambda | k_{m} | <1$. As $\|
  \tilde{K} \|_{L^{\infty}} \leqslant \sum_{m=1}^{\infty} | k_{m} |$, Theorem
  \ref{thm:201304021429} is proved. 

\begin{remark}
  For Onsager kernel we have $K ( \theta ) = | \sin   \theta |$, $k_{m} =-
  \frac{4}{\pi ( 4m^{2} -1 )}$. Theorems \ref{thm:201304021429} then gives
  $\lambda_{0} = \frac{\pi}{2}$. 
\end{remark}

\section{Bifurcation Analysis}\label{sec:201303261453}

We study the criticality and stability of bifurcated solutions from the
trivial solution. Our results generalize Proposition 4.4 and Corollary 4.5 in
{\cite{Lucia2010}}.

\begin{theorem}
  \label{thm:201502121119}Let $k_{m} <0$ satisfies $-k_{1} >-k_{2} > \cdots
  >0$. Then
  \begin{enumeratealpha}
    \item \textbf{(Bifurcation points)} two solutions bifurcate from the trivial solution at every
    $\lambda_{m} \assign - \frac{2}{k_{m}}$.
    
    \item \textbf{(Criticality)}
    \begin{itemize} 
      \item if $\frac{2k_{2m}}{k_{m}} <1$, both bifurcated solutions from
      $\lambda_{m}$ are supercritical;
      
      \item if $\frac{2k_{2m}}{k_{m}} >1$, both bifurcated solutions from
      $\lambda_{m}$ are subcritical. 
    \end{itemize}
    \item \textbf{(Stability)} the bifurcated solutions from $\lambda_{m} ,m \geqslant 2$ are
    unstable. The bifurcated solutions from $\lambda_{1}$ are stable if
    $\frac{2k_{2}}{k_{1}} <1$ and unstable if $\frac{2k_{2}}{k_{1}} >1$.
  \end{enumeratealpha}
  
\end{theorem}

\noindent\textbf{Proof.}
  
  \begin{enumeratealpha}
    \item At the trivial solution $V=0$ we have $\mathd \mu_{V} = \frac{\mathd
    \theta}{2 \pi}$. Following the calculation in Section
    \ref{subsec:201303271012}, the matrix $( a_{m n} )$ for the Jacobian $D 
    \Gamma$ is given by
    \begin{equation}
      a_{m n} = \left\{ \begin{array}{ll}
        - \lambda \frac{k_{m}}{2} & m=n\\
        0 & m \neq n
      \end{array} \right.
    \end{equation}
    and is thus diagonal.
    
    Denoting $\lambda_{m} =- \frac{2}{k_{m}}$ and $K=D  \Gamma$, we see that
    $\lambda_{m}$ is a simple characteristic value of $K$, and the dimensions
    of $\ker ( I- \lambda_{m} K )$ and $[ \tmop{Ran} ( I- \lambda_{m} K )
    ]^{\bot}$ are both $1$, which means $I- \lambda_{m} K$ is Fredholm with
    index zero.
    
    Furthermore, as $K=D  \Gamma ( 0 )$ and $\Gamma ( 0 ) =0$ we have
    \begin{equation}
      ( \Gamma -D  \Gamma ( 0 ) ) ( V ) =o ( \| V \| ) \hspace{2em} \text{ as
      } V \longrightarrow 0 \text{ in } H.
    \end{equation}
    Now setting $G ( \lambda ,V ) \assign - \lambda ( \Gamma -D  \Gamma ) ( V
    )$, we see that the problem we are solving, $V= \lambda \Gamma ( V )$,
    becomes $F ( \lambda ,V ) \assign V- \lambda K ( V ) +G ( \lambda ,V ) =0$
    with $G ( \lambda ,V ) =o ( \| V \| )$ as $V \longrightarrow 0$ uniformly
    in $\lambda$ near each $\lambda_{m}$. Therefore we can apply Theorem
    \ref{thm:201401141601} and Corollary \ref{cor:201402230858} (see Appendix
    \ref{subsec:nonana}) to conclude that two solutions bifurcate from the
    trivial solution at every $\lambda_{m} \assign - \frac{2}{k_{m}}$.
    
    \item Recall
    \begin{equation}
      \Gamma ( V ) = \frac{\int_{0}^{2 \pi} \tilde{K} ( \theta - \theta' )
      e^{-V ( \theta' )} \mathd \theta'}{\int_{0}^{2 \pi} e^{-V ( \theta )}
      \mathd \theta} .
    \end{equation}
    Expanding around the trivial solution we have
    \begin{equation}
      \Gamma ( V ) =T ( V ) +N ( V )
    \end{equation}
    where
    \begin{equation}
      T ( V ) =- \frac{1}{2 \pi}  \int_{0}^{2 \pi} \tilde{K} ( \theta -
      \theta' ) V ( \theta' ) \mathd \theta' ;
    \end{equation}
    \begin{equation}
      N ( V ) =- \frac{1}{2} T ( V^{2} ) + \frac{1}{6} T ( V^{3} ) - \left(
      \frac{1}{2 \pi}  \int_{0}^{2 \pi} \frac{V^{2}}{2} \right) T ( V ) +O (
      V^{4} ) ;
    \end{equation}
    Writing the orthonormal basis as
    \begin{equation}
      \phi_{n} ( \theta ) =c_{n} \cos ( 2n \theta ) ,
    \end{equation}
    where $c_{n} =1/ \sqrt{4n^{2} +1}$, we have
    \begin{equation}
      T ( \phi_{n} ) = \left( - \frac{k_{n}}{2} \right) \phi_{n} =
      \lambda_{n}^{-1} \phi_{n} .
    \end{equation}
    By Theorem \ref{thm:201401141601} and Corollary \ref{cor:201402230858} (see Appendix \ref{subsec:nonana}) we
    can write the bifurcated solution from $\lambda_{n}$ as
    \begin{equation}
      V=t \phi_{n} +t^{2} z
    \end{equation}
    where the $H^{1}$-inner product $\langle \phi_{n} ,z \rangle =0$.
    
    Next writing $\lambda = \lambda_{n} + \mu$ we have
    \begin{eqnarray}
      t^{2} z & = & \left( - \frac{2}{k_{n}} \right) t^{2} T ( z ) + \mu t^{2}
      T ( z ) + \mu t \left( - \frac{k_{n}}{2} \right) \phi_{n} \nonumber\\
      &  & + ( \lambda_{n} + \mu )  \left[ - \frac{t^{2}}{4} 
      \frac{c_{n}^{2}}{c_{2n}}  \left( - \frac{k_{2n}}{2} \right) \phi_{2n}
      -t^{3} T ( \phi_{n} z ) \right] \nonumber\\
      &  & + ( \lambda_{n} + \mu )  \left[ \frac{t^{3}}{24} 
      \frac{c_{n}^{3}}{c_{3n}}  \left( - \frac{k_{3n}}{2} \right) \phi_{3n} -
      \frac{t^{3}}{8} c_{n}^{2}  \left( - \frac{k_{n}}{2} \right) \phi_{n}
      \right] +O ( t^{4} ) .  \label{eq:201401231113}
    \end{eqnarray}
    Taking $H^{1}$-inner product with $\phi_{n}$ and using the facts that
    $\langle \phi_{n} ,z \rangle =0$, $\langle T ( z ) , \phi_{n} \rangle =0$
    we reach
    \begin{equation}
      \mu =t^{2} \lambda_{n}  \left( - \frac{k_{n}}{2} \right)^{-1}  \left[
      \langle T ( \phi_{n} z ) , \phi_{n} \rangle + \left( - \frac{k_{n}}{2}
      \right)  \frac{c_{n}^{2}}{8} \right] +O ( t^{3} ) .
    \end{equation}
    Consequently $\mu =O ( t^{2} )$ and (\ref{eq:201401231113}) can be
    simplified to
    \begin{eqnarray}
      t^{2} z & = & \left( - \frac{2}{k_{n}} \right) t^{2} T ( z ) + \mu t
      \left( - \frac{k_{n}}{2} \right) \phi_{n} \nonumber\\
      &  & + \left( - \frac{2}{k_{n}} \right)  \left[ - \frac{t^{2}}{4} 
      \frac{c_{n}^{2}}{c_{2n}}  \left( - \frac{k_{2n}}{2} \right) \phi_{2n}
      -t^{3} T ( \phi_{n} z ) \right] \nonumber\\
      &  & + \left( - \frac{2}{k_{n}} \right)  \left[ \frac{t^{3}}{24} 
      \frac{c_{n}^{3}}{c_{3n}}  \left( - \frac{k_{3n}}{2} \right) \phi_{3n} -
      \frac{t^{3}}{8} c_{n}^{2}  \left( - \frac{k_{n}}{2} \right) \phi_{n}
      \right] +O ( t^{4} ) .  \label{eq:201401231116}
    \end{eqnarray}
    To obtain the sign of $\mu$, we need to calculate $\langle T ( \phi_{n} z
    ) , \phi_{n} \rangle$. Writing
    \begin{equation}
      z= \sum_{k=2}^{\infty} z_{k} \phi_{k} .
    \end{equation}
    We have
    \begin{equation}
      \langle T ( \phi_{n} z ) , \phi_{n} \rangle = \left( - \frac{k_{n}}{2}
      \right) z_{2n}  \frac{c_{2n}}{2} .
    \end{equation}
    Now we calculate $z_{2n}$. Taking $H^{1}$-inner product of
    (\ref{eq:201401231116}) with $\phi_{2n}$, we finally reach
    \begin{equation}
      z_{2n} = \frac{\gamma_{n}}{\gamma_{n} -1}  \frac{c_{n}^{2}}{4c_{2n}}
    \end{equation}
    where $\gamma_{n} \assign \frac{k_{2n}}{k_{n}}$.
    
    Putting things together, we have
    \begin{equation}
      \mu =t^{2} \lambda_{n}  \left[ z_{2n}  \frac{c_{2n}}{2} +
      \frac{c_{n}^{2}}{8} \right] +O ( t^{3} ) =t^{2} \lambda_{n} 
      \frac{c_{n}^{2}}{8}  \frac{2 \gamma -1}{\gamma -1} +O ( t^{3} )
    \end{equation}
    and consequently the bifurcation is super-critical if $( 2 \gamma_{n} -1 )
    / ( \gamma_{n} -1 ) >0$ and sub-critical if $( 2 \gamma_{n} -1 ) / (
    \gamma_{n} -1 ) <0$. The conclusion of the theorem thus follows. \
    
    \item It is clear that the trivial solution is stable for $\lambda <
    \lambda_{1}$ and unstable for $\lambda > \lambda_{1}$. Therefore the
    bifurcated solutions from $\lambda_{m}$ with $m \geqslant 2$ are unstable,
    independent of their criticality.
    
    For the bifurcations from $\lambda_{1}$, we check that the assumptions of
    Theorem \ref{thm:201402111013} (see Appendix \ref{subsec:nonana}) are all
    satisfied at $\lambda_{1}$.
    
    First recall that a bifurcation point $\mu_{0}$ is ``regular'' if the
    linearized operator is invertible ``for all $\mu$ sufficiently close to
    $\mu_{0}$ but $\mu \neq \mu_{0}$'' ({\cite{Sattinger1971}}). As the
    linearized operator $D  \Gamma$ is diagonal, we easily see that all the
    bifurcation points under discussion are regular.
    
    Next we check the smoothness conditions for Theorem
    \ref{thm:201402111013}. The nonlinear remainder term $N ( \lambda ,V )$ is
    given by
    \begin{eqnarray}
      N ( \lambda ,V ) & = & \lambda ( D  \Gamma ( 0 ) - \Gamma ) ( V )
      \nonumber\\
      & = & - \lambda \left[ \frac{1}{2 \pi}  \int_{0}^{2 \pi} \tilde{K} (
      \theta - \theta' ) V ( \theta' ) \mathd \theta' + \frac{\int_{0}^{2 \pi}
      \tilde{K} ( \theta - \theta' ) e^{-V ( \theta' )} \mathd
      \theta'}{\int_{0}^{2 \pi} e^{-V ( \theta )} \mathd \theta} \right] . 
    \end{eqnarray}
    It is easy to see that $N ( \lambda ,V )$ is Fr{\'e}chet differentiable as
    \begin{equation}
      \delta V \rightarrow 0 \text{ in } H \Longrightarrow \delta V
      \rightarrow 0 \text{ in } L^{\infty} \Longrightarrow e^{- ( V+ \delta V
      )} \longrightarrow e^{-V} \text{ uniformly}
    \end{equation}
    thanks to the embedding $H \longhookrightarrow L^{\infty} ( [ 0,2 \pi ]
    )$. Similarly we can show that $N ( \lambda ,V )$ is twice Fr{\'e}chet
    differentiable.
    
    Finally we define $N_{1} ( \lambda ,V, \alpha ) \assign \alpha^{-2} N (
    \lambda , \alpha V )$ and prove that it is Fr{\'e}chet differentiable in
    $\lambda ,V$ and $\alpha$. It is obvious that $N_{1}$ is Fr{\'e}chet
    differentiable in $\lambda$ and $V$. To see that it is also Fr{\'e}chet
    differentiable in $\alpha$, we write
    \begin{equation}
      N_{1} ( \lambda ,V, \alpha ) =- \lambda \int_{0}^{2 \pi} \tilde{K} (
      \theta - \theta' ) R ( \theta' , \alpha ) \mathd \theta'
    \end{equation}
    where
    \begin{equation}
      R ( \theta , \alpha ) \assign \alpha^{-2}  \left[ \frac{-1+ \alpha V (
      \theta )}{2 \pi} + \frac{e^{- \alpha V ( \theta )}}{\int_{0}^{2 \pi}
      e^{- \alpha V ( \theta )} \mathd \theta} \right] .
      \label{eq:201502181352}
    \end{equation}
    Using Taylor expansion with integral form of remainder, we have
    \begin{eqnarray}
      e^{- \alpha V ( \theta )} & = & 1- \alpha V ( \theta ) +
      \int_{0}^{\alpha V ( \theta )} ( \alpha V ( \theta ) -t ) e^{-t} \mathd
      t \nonumber\\
      & = & 1- \alpha V ( \theta ) + \alpha^{2}  \int_{0}^{V ( \theta )} ( V
      ( \theta ) -s ) e^{- \alpha s} \mathd s. 
    \end{eqnarray}
    Substituting this into (\ref{eq:201502181352}) we have
    \begin{equation}
      R ( \theta , \alpha ) = \frac{( -1+ \alpha V ( \theta ) )  \int_{0}^{2
      \pi} M ( \theta , \alpha ) \mathd \theta +2 \pi M ( \theta , \alpha )}{2
      \pi \int_{0}^{2 \pi} e^{- \alpha V ( \theta )} \mathd \theta}
    \end{equation}
    where $M ( \theta , \alpha ) \assign \int_{0}^{V ( \theta )} ( V ( \theta
    ) -s ) e^{- \alpha s} \mathd s$. We see that clearly $R ( \theta , \alpha
    )$ is differentiable in $\alpha$ and so is $N_{1} ( \lambda ,V, \alpha )$.
    
    Thus we can apply Theorem \ref{thm:201402111013} to immediately conclude:
    \begin{itemize}
      \item If $\frac{2k_{2}}{k_{1}} <1$, then the bifurcated solutions from
      $\lambda_{1}$ are stable;
      
      \item If $\frac{2k_{2}}{k_{1}} >1$, then the bifurcated solutions from
      $\lambda_{1}$ are unstable. 
    \end{itemize}
  \end{enumeratealpha}
\vspace{3mm}

\begin{remark}
  All the bifurcations from the trivial solution for the Onsager model are
  super-critical. For Onsager kernel the bifurcation values are $\lambda_{m} =
  \frac{( 4m^{2} -1 ) \pi}{2}$, $m=1,2,3, \ldots$. We see that the first
  bifurcation value is $\frac{3 \pi}{2}$. Thus there is a gap between it and
  the uniqueness region $\lambda < \lambda_{0} = \frac{\pi}{2}$.
\end{remark}

\appendix\section{Auxiliary Lemmas and Known Theorems}\label{app:201303261454}

\subsection{Classical Results from Nonlinear Analysis}\label{subsec:nonana}

The following classical results in nonlinear analysis are crucial in our
analysis.

\begin{theorem}[{\cite{Deimling1985}}, Theorem 28.3]
  \label{thm:201401141601}Let $X$ be a real Banach space, $K \in L ( X )$,
  $\Omega \subset \mathbbm{R} \times X$ a neighborhood of $( \lambda_{0} ;0 )$
  and $G: \Omega \mapsto X$ such that $G_{\lambda} ,G_{x} ,G_{\lambda  x}$ are
  continuous on $\Omega$. Suppose also that
  \begin{enumeratealpha}
    \item $G ( \lambda ,x ) =o ( \| x \| )$ as $x \longrightarrow 0$ uniformly
    in $\lambda$ near $\lambda_{0}$.
    
    \item $I- \lambda_{0} K$ is Fredholm of index zero and $\lambda_{0}$ is a
    simple characteristic value of $K$.
  \end{enumeratealpha}
  Then $( \lambda_{0} ;0 )$ is a bifurcation point for $F ( \lambda ,x ) =x-
  \lambda K+G ( \lambda ,x ) =0$ and there is a neighborhood $U$ of $(
  \lambda_{0} ;0 )$ such that
  \begin{equation}
    F^{-1} ( 0 ) \cap U= \{ ( \lambda_{0} + \mu ( t ) ,tv+tz ( t ) ) : | t | <
    \delta \} \cup \{ ( \lambda ;0 ) : ( \lambda ;0 ) \in U \}
  \end{equation}
  for some $\delta >0$, with continuous functions $\mu ( \cdot )$ and $z (
  \cdot )$ such that $\mu ( 0 ) =0,z ( 0 ) =0$ and the range of $z ( \cdot )$
  is contained in a complement of $N ( I- \lambda_{0} K ) = \tmop{span} \{ v
  \}$.
\end{theorem}

\begin{corollary}[{\cite{Deimling1985}}, Corollary 28.1]
  \label{cor:201402230858}Let the hypotheses of Theorem 28.3 be fulfilled. If
  $G$ is $C^{k}$ near $( \lambda_{0} ;0 )$ for some $k \geqslant 2$ then the
  functions $\mu ( \cdot ) ,z ( \cdot ) \nocomma$, defining the branches of
  nontrivial zeros, are $C^{k-1}$. If $G$ is real (or complex) analytic then
  $\mu ( \cdot )$ and $z ( \cdot )$ are real (or complex) analytic.
\end{corollary}

We also made use of the following result by Sattinger.

\begin{theorem}[{\cite{Sattinger1971}}, Theorem 4.2]
  \label{thm:201402111013}Let $( \mu_{0} ,0 )$ be a regular bifurcation point
  of (3.1) and let $N$ be twice continuously Fr{\'e}chet differentiable, with
  $N ( \mu , \alpha u ) = \alpha^{2} N_{1} ( \mu ;u; \alpha )$ where $N_{1}$
  is Fr{\'e}chet differentiable in $\mu ,u$ and $\alpha$. Then the
  supercritical bifurcating solutions are stable and subcritical bifurcating
  solutions are unstable. 
\end{theorem}

\subsection{Properties of $H$ and $\Gamma$}\label{app:201303270947}

\noindent{\textbf{Proof of Lemma \ref{lem:201304161059}.} We first prove that $H$, as defined in (\ref{eq:201303140924}),
\begin{equation}
  H \assign \left\{ V ( \theta ) \in H^{1} ( [ 0,2 \pi ] ) ;V ( \theta ) =V (
  \theta + \pi )  a.e.; \int_{0}^{2 \pi} V ( \theta ) \mathd \theta =0;V (
  \theta ) =V ( 2 \pi - \theta )  a.e. \right\} . \label{eq:201502111431}
\end{equation}
is a Hilbert space.

  Since $H$ is a subspace of the Hilbert space $H^{1} ( [ 0,2 \pi ] )$, all we
  need to show is that it is closed in the topology of $H^{1}$, which is
  trivial.
  
  Next it is easy to check that $\Gamma ( V )$ satisfies the 2nd, 3rd, and
  4th requirements in (\ref{eq:201502111431}). To show that it is in $H^{1}$,
  we calculate
  \begin{eqnarray}
    \| \Gamma ( V ) \|_{L^{2}}^{2} = \int_{0}^{2 \pi} [ \Gamma ( V ) ( \theta
    ) ]^{2} \mathd \theta & = & \int_{0}^{2 \pi}  \left[ \frac{\int_{0}^{2
    \pi} \tilde{K} ( \theta - \theta' ) e^{-V ( \theta' )} \mathd
    \theta'}{\int_{0}^{2 \pi} e^{-V ( \theta )} \mathd \theta} \right]^{2}
    \mathd \theta \nonumber\\
    & \leqslant & 2 \pi \| \tilde{K} \|_{L^{\infty}}^{2} =2 \pi \| K- \bar{K}
    \|_{L^{\infty}}^{2} < \infty . \label{eq:201303261508} 
  \end{eqnarray}
  Similarly, we have $\left\| \frac{\mathd}{\mathd \theta} \Gamma ( V )
  \right\|_{L^{2}}^{2} \leqslant 2 \pi \| K' \|_{L^{\infty}}^{2} < \infty$. 
Thus ends the proof of Lemma \ref{lem:201304161059}. 
\vspace{3mm}

\noindent\textbf{Proof of Lemma \ref{lem:201303270911}.} Next we prove the continuity and compactness of the operator $\Gamma$. Recall
that $\Gamma$ is defined in (\ref{eq:201502111432}) as
\begin{equation}
  \Gamma ( V ) ( \theta ) = \frac{\int_{0}^{2 \pi} \tilde{K} ( \theta -
  \theta' ) e^{-V ( \theta' )} \mathd \theta'}{\int_{0}^{2 \pi} e^{-V ( \theta
  )} \mathd \theta} . \label{eq:201304021416}
\end{equation}

  \begin{itemize}
    \item Continuity.
    
    Let $\delta V \longrightarrow 0$ in $H$. We first show that $\Gamma ( V+
    \delta V ) ( \theta ) \longrightarrow \Gamma ( V ) ( \theta )$ in $L^{2}$.
    Thanks to the embedding $H^{1} ( [ 0,2 \pi ] ) \longhookrightarrow
    L^{\infty} ( [ 0,2 \pi ] )$, we have $\delta V \longrightarrow 0$ in
    $L^{\infty}$. Consequently $e^{- ( V+ \delta  V )} \longrightarrow e^{-V}$
    uniformly and it follows that $\Gamma ( V+ \delta V ) ( \theta )
    \longrightarrow \Gamma ( V ) ( \theta )$ uniformly and the conclusion
    follows.
    
    Next we show that $\frac{\mathd}{\mathd \theta} \Gamma ( V+ \delta V ) (
    \theta ) \longrightarrow \frac{\mathd}{\mathd \theta} \Gamma ( V ) (
    \theta )$ in $L^{2}$. We calculate
    \begin{equation}
      \frac{\mathd}{\mathd \theta} \Gamma ( V ) ( \theta ) = \frac{-
      \int_{0}^{2 \pi} \tilde{K} ( \theta - \theta' ) e^{-V ( \theta' )} V' (
      \theta' ) \mathd \theta'}{\int_{0}^{2 \pi} e^{-V ( \theta )} \mathd
      \theta} .
    \end{equation}
    As $e^{- ( V+ \delta V )} \longrightarrow e^{-V}$ uniformly, $e^{- ( V+
    \delta V )}  ( V' + \delta V' ) \longrightarrow e^{-V} V'$ in $L^{2}$
    which together with $\tilde{K} \in L^{\infty}$ implies
    $\frac{\mathd}{\mathd \theta} \Gamma ( V+ \delta V ) ( \theta )
    \longrightarrow \frac{\mathd}{\mathd \theta} \Gamma ( V ) ( \theta )$ in
    $L^{\infty}$ and consequently also in $L^{2}$.
    
    \item Compactness. Assume $\Omega \subset B_{R}$ be a bounded closed
    subset of $H$, where $B_{R}$ denotes the ball with radius $R$ in $H$. It
    suffices to show that there are operators $\Gamma_{n} \longrightarrow
    \Gamma$ whose ranges are finite dimensional (see e.g.
    {\cite{Nirenberg2001}}). By the Weierstrass approximation theorem, for any
    $n \in \mathbbm{N}$, there is $K_{n} ( \theta ) = \sum_{i=1}^{l_{n}} [
    a_{n i} \cos ( m_{n i} \theta ) +b_{n i} \sin ( m_{n i} \theta ) ]$ such
    that
    \begin{equation}
      | \tilde{K} ( \theta ) -K_{n} ( \theta ) | < \frac{1}{n} \hspace{2em}
      \forall \theta \in [ 0,2 \pi ] .
    \end{equation}
    Now we define
    \begin{equation}
      \Gamma_{n} ( V ) ( \theta ) \assign \frac{\int_{0}^{2 \pi} K_{n} (
      \theta - \theta' ) e^{-V ( \theta' )} \mathd \theta'}{\int_{0}^{2 \pi}
      e^{-V ( \theta )} \mathd \theta}
    \end{equation}
    It is easy to check that
    \begin{equation}
      \Gamma_{n} ( V ) ( \theta ) \in \tmop{span} \{ \cos ( m_{n i} \theta ) ,
      \sin ( m_{n i} \theta ) \}_{i=1}^{l_{n}}
    \end{equation}
    for any $V ( \theta ) \in H$ which means the range of $\Gamma_{n}$ is
    finite dimensional.
    
    Finally check
    \begin{eqnarray}
      \| \Gamma_{n} - \Gamma \|_{\Omega \mapsto H} & \leqslant & \sup_{\| V
      \|_{H^{1}} \leqslant R} \left\| \frac{\int_{0}^{2 \pi} ( \tilde{K} (
      \theta - \theta' ) -K_{n} ( \theta - \theta' ) ) e^{-V ( \theta' )}
      \mathd \theta'}{\int_{0}^{2 \pi} e^{-V ( \theta )} \mathd \theta}
      \right\|_{H^{1}} \nonumber\\
      & \leqslant & \sup_{\| V \|_{L^{2}} \leqslant R} \left\|
      \frac{\int_{0}^{2 \pi} ( \tilde{K} ( \theta - \theta' ) -K_{n} ( \theta
      - \theta' ) ) e^{-V ( \theta' )} \mathd \theta'}{\int_{0}^{2 \pi} e^{-V
      ( \theta )} \mathd \theta} \right\|_{L^{2}} \nonumber\\
      &  & + \sup_{\| V \|_{H^{1} \leqslant R}} \left\| \frac{\int_{0}^{2
      \pi} ( \tilde{K} ( \theta - \theta' ) -K_{n} ( \theta - \theta' ) )
      e^{-V ( \theta' )} V' ( \theta' ) \mathd \theta'}{\int_{0}^{2 \pi} e^{-V
      ( \theta )} \mathd \theta} \right\|_{L^{2}} \nonumber\\
      & \leqslant & C ( R ) [ \| \tilde{K} -K_{n} \|_{L^{\infty}} + \|
      \tilde{K} -K_{n} \|_{L^{2}} ] \leqslant \frac{C ( R )}{n} . 
    \end{eqnarray}
    The calculation is similar to that in the proof of Lemma
    \ref{lem:201304161059} and is omitted here. The arbitrariness of $n$ now
    gives the desired result. 
  \end{itemize}

\subsection{A Gr{\"u}ss type inequality}\label{subsec:201502101458}

The following Gr{\"u}ss type inequality will play a crucial role in the
proofs.

\begin{lemma}
  \label{lem:201304181013}Let $\mu$ be a probability measure over a domain
  $\Omega$. Let $f,g \in L^{\infty} ( \Omega )$ satisfy $a \leqslant f
  \leqslant A,b \leqslant g \leqslant B$. Then
  \begin{equation}
 \lvert
\int_{\Omega} f ( x ) g ( x ) \mathd \mu - \left( \int_{\Omega} f (x) \mathd \mu \right)  
\left(\int_{\Omega} g ( x ) \mathd \mu \right) 
\rvert
\leqslant \frac{( A-a )  ( B-b )}{4} .
  \end{equation}
  in particular we have
  \begin{equation}
    \lvert \int_{\Omega} f ( x ) g ( x ) \mathd \mu - \left( \int_{\Omega} f (
    x ) \mathd \mu \right)  \left( \int_{\Omega} g ( x ) \mathd \mu \right)
    \rvert \leqslant \| f \|_{L^{\infty}}  \| g \|_{L^{\infty}} .
  \end{equation}
\end{lemma}

\noindent\textbf{Proof.}
  This is a simple generalization of the classical Gr{\"u}ss inequality. The
  proof is almost identical to that in {\cite{Dragomir2000}} and is therefore
  omitted.

\paragraph{Acknowledgment}A. Niksirat and X. Yu are partially supported by a
grant from NSERC. The authors would like to thank Prof. Pingwen Zhang for
suggestion of the problem and valuable discussions. The authors would also
like to thank the anonymous referee and the associate editor for valuable
comments and advice that greatly improve the organization and presentation of
the paper.

\end{document}